\documentclass[smallextended,envcountsect]{svjour3}
\smartqed
\usepackage{graphicx}
\usepackage{amsfonts}
\usepackage{amsmath}
\usepackage[usenames,dvipsnames]{pstricks}
\usepackage{epsfig}
\usepackage{bbding}

\newcommand{\bc}{\begin{center}}
\newcommand{\br}{\begin{right}}
\newcommand{\ec}{\end{center}}
\newcommand{\be}{\begin{equation}}
\newcommand{\ee}{\end{equation}}

\newcommand{\pr}{\parallel}
\newcommand{\vl}{\mid}

\newcommand{\rar}{\rightarrow}
\newcommand{\grad}{\nabla}
\newcommand{\lapl}{\triangle}

\newcommand{\p}{\partial}

\newcommand{\bint}{\mbox{int} \,}

\newcommand{\meq}{\geq}

\newcommand{\al}{\alpha}

\newcommand{\eps}{\varepsilon}

\newtheorem{thm}{Theorem}[section]

\newtheorem{cor}{Corollary}[section]

\newtheorem{rem}{Remark}[section]

\newcommand{\incl}{\subset}

\journalname{JOTA}

\begin{document}

\title{The accelerated method for finding the minimum of a
nonsmooth finite convex function}


\titlerunning{The accelerated method of finding minimum of a
convex  function}

\author{Igor M. Prudnikov}

\institute{Igor Mihailovich Prudnikov \at
             Scientific Center of Smolensk Federal Medical  University, Smolensk, Russia, 214000\\
              \email{pim\underline{ }10@hotmail.com}
}

\date{Received: date / Accepted: date}

\maketitle

\begin{abstract}
The goal of the paper is development of an  optimization method
with a superlinear  convergence rate for a nonsmooth convex
function. For optimization an approximation is used that is
similar to the Steklov integral averaging. The difference is that
averaging is performed over a variable-dependent set, that is
called a set-valued mapping (SVM) satisfying simple conditions.
The novelty of the approach is that with such an approximation we
obtain twice continuously differentiable convex functions, for the
optimization of which  of the second order methods are used. The
rate of convergence of the method is estimated.
\end{abstract}



\section{Introduction}

Nonsmooth (non-differentiable) or insufficiently smooth functions
are widely used in economics, data processing, control theory,
artificial intelligence and other fields. An example of such
functions are functions obtained by performing operations minimum
or maximum.

Nonsmooth functions may not have derivatives at some points.  It
is known that the Lipschitz function is differentiable almost
everywhere (a.e.) in $R^n$ \cite{rademacher}.  Generalized
gradients are used instead of  gradients at the points of
non-differentiability of a function. The optimization methods of
these functions are different from the optimization methods  of
smooth (differentiable) functions.

In this paper the author continues research related to the
construction of an optimization method of Lipschitz functions
using the Steklov integrals and similar integrals, when a set,
over which averaging is taken, is a function of a variable.

This approach gives twice differentiable functions, whose
stationary points coincide with the stationary points of the
original function in contrast to the case when averaging is doing
over sets independent of $x$. For such functions second-order
optimization methods can be used that are tested for arbitrary
convex functions with an estimate of the convergence rate.

If we have discontinuous gradients as functions of variables, then
it is very difficult to construct optimization methods and
estimate their convergence rates in the general case. Using the
polynomial approximation of an original function and transition to
optimization of a smooth function by the known methods
\cite{pshenichnyidanilin} does not allow to solve the optimization
problem, since this way leads to appearance of new extremum points
located far from the extremum points of the original function.

Separation of fictitious extremum points from real ones is as
complex a problem as the initial one. Therefore, the development
of the theory of nonsmooth functions went along the path of
developing its own methods, based on the properties of generalized
gradient of  Lipschitz functions. Here it is worth mentioning the
articles  \cite{pshenichnyidanilin} - \cite{clarke} N. Z. Shor, B.
N. Pshenichny, V. F. Demyanova, E.A. Nurminsky, F. Clark, R.T.
Rokafellar, L.N. Polyakova.

To construct accelerated optimization methods for nonsmooth
functions, it is necessary to determine the constructions to which
second-order optimization methods are applicable. But to perform
the latter, it is necessary to determine such constructions for
which the extremum points do not disappear, and new ones do not
appear.

The paper proposes exactly this method for smoothing of nonsmooth
functions. The resulting function will be continuously
differentiable. If we apply again the averaging operation to it,
then we will have a twice differentiable function.

If we apply averaging over sets depending on the variable $ x $,
then we obtain a continuously differentiable function whose
stationary points coincide with the stationary points of the
original function. If we repeat the averaging procedure, we get
twice differentiable functions to which second-order optimization
methods with accelerated convergence can be applied.

It is possible to move with the help of the defined functions from
the local optimization of non-smooth functions to local
optimization of smooth functions, and also to estimate the rate of
convergence to an extremum point, that is definitely  important,
because it is possible to develop accelerated optimization methods
for functions with discontinuous gradients. Similar constructions
as far as the author knows, nobody has proposed previously.

\section{Smoothing integral functions}

Let $f(\cdot): R^n \rar R$  be a Lipschitz function with a
constant $L$, $x_* $ is its local minimum (maximum) in $R^n$. As
it is known, the necessary extremum condition at the point $x_*$
for the Lipschitz function $f(\cdot)$ is zero belongs to the
Clarke subdifferential $\p_{CL} f(\cdot),$  calculated at the
point $x_* $, i.e.
$$
0 \in \p_{CL} f(x_*).
$$
Any point, for which this condition is correct, is called a
stationary point. Not all stationary points are minimum or maximum
points.

Let us take an arbitrary convex compact set $D \incl R^n$, $0 \in
\bint D$. We introduce the definition of the $\eps(D)$ stationary
point.

{\bf Definition 1.} A point $x_{\eps}$ is called  the $\eps(D)$
stationary point of the function $f(\cdot)$ , if the set $
x_{\eps} + D$ includes a stationary point of the function
$f(\cdot)$.

This definition agrees with the definition of the $\eps$
stationary point for the convex functions \cite{rocafellar},
because for the strongly convex functions the distance from the
$\eps$ stationary point to the minimum can be evaluated by
difference of values of the function $f(\cdot)$ calculated at
these points.

Define the function $\varphi(\cdot): R^n \rar R$ \be \varphi(x)=
\frac{1}{\mu(D)} \int_D f(x+y)d y, \label{newopt1a} \ee where
$\mu(D)$ is the measure of the domain $D$, $\mu(D) >0$.

Obviously, $\varphi(\cdot)$  is continuous. Let us show that
$\varphi(\cdot)$ is a Lipschitz function with the Lipschitz
constant equaled to the Lipschitz constant of the function
$f(\cdot)$. Really,
$$
\vl \varphi(x_1) - \varphi(x_2) \vl \leq \frac{1}{\mu(D)} \int_D
\vl f(x_1+y)-f(x_2+y) \vl d y  \leq \frac{1}{\mu(D)} \int_D L \|
x_1 -x_2 \| dy \leq
$$
$$
\leq L \| x_1-x_2 \| \,\,\,\,\,\, x_1,x_2 \in R^n.
$$
The function $f(\cdot)$  is Lipschitz, and therefore it is
differentiable a.e. in $R^n$ \cite{rademacher}. Let $N(f)$ denote
the set of points of differentiability of the function $f(\cdot)$
in $R^n$. It is known that $N(f)$ is everywhere dense in $R^n$
and, in particular, in $D$, because of $\mu(D)>0$ by assumption.

The following theorem was proved in \cite{proudintegapp1}.

\begin{thm} For an arbitrary Lipschitz function $f(\cdot): R^n \rar
R$ the function
$$
\varphi(x)= \frac{1}{\mu(D)} \int_D f(x+y)d y,
$$
where $D$ is any domain in $R^n, 0 \in int D, \mu(D)$ is the
measure of the domain $D$, $\mu(D) >0$,  is a continuously
differentiable function with the derivative
$$
\varphi'(x)=\frac{1}{\mu(D)} \int_D f'(z+x )d z.
$$
\end{thm}

\begin{rem} We use here the Lebesque integration. \end{rem}
\begin{rem} The derivatives of the function $f(\cdot)$ are taken at
those points where they exist.\end{rem}

It was also proved in \cite{proudintegapp1} that if $f(\cdot)$ is
Lipschitz, then $\varphi'(\cdot)$ is also Lipschitz function.

Consider the function
$$
\phi(x)= \frac{1}{\mu(D)} \int_D \varphi(x+y)d y.
$$
Since $\varphi(\cdot)$ is Lipschitz, we will have \be \phi'(x)=
\frac{1}{\mu(D)} \int_D \varphi'(z+x )d z. \label{newopt4} \ee

Since $\varphi'(\cdot)$ is continuous, $\phi(\cdot)$ is a
continuously differentiable function. As soon as $\varphi'(\cdot)$
is Lipschitz, we can differentiate (\ref{newopt4}). As a result,
we will have \be \phi''(x)= \frac{1}{\mu(D)} \int_D \varphi''(z+x
)dz, \label{newopt5} \ee i.e. $\phi(\cdot)$ is a twice
continuously differentiable function.

It can be shown \cite{proudintegapp} that the function
$\phi''(\cdot)$ is Lipschitz with a constant $\tilde{L}$,
depending on the set $D$. If $D$ is a ball or a cube in
$\mathbb{R}^n $, then  we can take $\tilde{L}=\frac{2L}{d^2}$,
where $d$ is the diameter of the set $D$, $L$ is the Lipschitz
constant of the function $f(\cdot)$.

\begin{rem} The integration in (\ref{newopt5}) is understood, as before, in
the sense of Lebesgue.
\end{rem}

If $x$ is a point of local maximum or minimum of the function
$f(\cdot)$, then for sufficiently small $r>0$ and
$D=S^{n-1}_r(0)=\{ z \in R^n \mid \| z \| \leq r \}$ the point $x
$ is also the local minimum or maximum point of the function
$\varphi(\cdot)$. But unlike the function $f(\cdot)$ the function
$\varphi(\cdot)$ is continuously differentiable. Similar thing is
true for the function $\phi(\cdot)$,  i.e. the point $x$ is a
point of local minimum or maximum of the function $\phi(\cdot)$.
But unlike the functions $f(\cdot)$ and $\varphi(\cdot)$ the
function $\phi(\cdot)$  is twice continuously differentiable,
matrix of the second mixed derivatives of which satisfies to the
Lipschitz condition. To optimize $\phi(\cdot)$ we can use the
methods of second order.

The functions $\varphi(\cdot)$ and  $\phi(\cdot)$ also retain many
properties of the function $f(\cdot)$. An important property for
applications of the functions $\varphi(\cdot)$ and $\phi(\cdot)$
is that if $f(\cdot)-$ is convex with respect to all or some
variables, then $\varphi(\cdot)$ and $\phi(\cdot)$ are also convex
with respect to the same variables \cite{proudintegapp}.

Let us see which stationary points the function $\varphi(\cdot)$
has. According to the formula (\ref{newopt4}), the stationary
point $x_*$ of the function $\varphi(\cdot)$ is such a point, for
which \be \varphi'(x_*)= \frac{1}{\mu(D)} \int_D f'(z+ x_* )d z
=0. \label{newopt6} \ee We will show that  the stationary point of
the function $f(\cdot)$ belongs to the set $x_*+D$.

The integral in (\ref{newopt6}) can be represented with any degree
of accuracy $\delta >0$ in the form \be \frac{1}{\mu(D)}
\sum_{i=1}^N f'(z_i+x_*) \mu (D_i), \label{newopt7} \ee where
$N=N(\delta)$, $D_i \incl D, i \in 1:N, $ are subregions of the
set $D$, $\mu (D_i)$ are their measures,
$$
\sum_{i=1}^N \mu (D_i) = \mu (D).
$$
The sum (\ref{newopt7}) is the convex hull of the vectors
$f'(z_i+x_*)$. Really, \be \frac{1}{\mu(D)} \sum_{i=1}^N
f'(z_i+x_*) \mu (D_i) = \sum_{i=1}^N \frac{\mu (D_i)}{\mu(D)}
f'(z_i+x_*)= \sum_{i=1}^N  \al_i f'(z_i+x_*), \label{newopt8} \ee
where $\alpha_i = \frac{\mu (D_i)}{\mu(D)}, \alpha_i \geq 0,$ and
$\sum_{i=1}^N \al_i = 1.$

According to the equality (\ref{newopt6}), the sum (\ref{newopt8})
can be made arbitrarily small for large $N=N(\delta)$ (for small
$\delta$ ). Since the convex hull of any vectors is a closed set
and the convex hull of generalized gradients is a collinear vector
to some generalized gradient of the function $f(\cdot)$ at a point
$x_* + \bar{z} \in x_*+D $, $ \bar{z} \in D$, we obtain that the
sum (\ref{newopt8}) is a vector tending to zero generalized
gradient as $N \rar \infty$. In other words, there exists a point
$x_* + \bar{z} \in x_*+D $, with a zero generalized gradient of
the function $f(\cdot)$.

Therefore, the stationary point $ x_* + \bar{z}$ of the function
$f(\cdot)$ belongs to the set $ x_* + D$. Hence, by definition,
$x_*$ is a $\eps(D)$ stationary point. Thus, the following theorem
is proved.

\begin{thm} All stationary points of the function $\varphi(\cdot)$
are the $\eps(D)$ stationary points of the function $f(\cdot)$.
\end{thm}

Similar reasoning is true for the function $\phi(\cdot)$.

\begin{cor} All stationary points of the function $\phi(\cdot)$
are the $\eps(D)$ stationary points of the function
$\varphi(\cdot)$ or the $\eps(2D)$stationary points of the
function $f(\cdot)$.
\end{cor}

\begin{cor}  If $x_*$ is a local minimum point of the function $f(\cdot)$,
for which there exists a neighborhood $S, x_* \in \bint S,$ where
$$
f(z) \meq f(x_*) \,\,\, \forall z \in S,
$$
then there exists a convex compact set $D$ and a point $y \in S$,
where $\varphi'(y)=0$ and $x_* \in y+D \incl S$, i.e. the point
$y$ is the $\eps(D)$ stationary point of the function $f(\cdot)$.
\end{cor}

The same is true for the local maximum point of the function
$f(\cdot)$.

To find the $ \eps (2D) $ stationary points of the function $ f
(\cdot) $, we must apply second-order optimization methods for the
function $ \phi (\cdot) $. A numerical optimization method will be
given with the rate of convergence to a stationary point of the
function $ f (\cdot) $ faster than any geometric progression.


\section{Search algorithm for stationary points of the Lipschitz function $f(\cdot)$}


Let us take a sequence of sets $\{ D_s \}, s=1,2,\dots $ with
non-empty interior  whose diameters $d(D_s)$ tends to zero with
increasing $s$. Let be $D_s=B^n_{r_s}(0)=\{ v \in \mathbb{R}^n \vl
\| v \| \leq r_s \}$ for $r_s \rar +0$ as $s \rar \infty$. We
introduce a sequence of the functions
$$
\varphi_s(x)= \frac{1}{\mu(D_s)} \int_{D_s} f(x+y)d y
$$
and
$$
\Phi_s(x)= \frac{1}{\mu(D_s)} \int_{D_s} \varphi_s (x+y)d y.
$$
Let the inequality $ \| \Phi''_s (\cdot) \| \leq {L}_{s}$ be true
for the matrix of the second mixed derivatives of the function
$\Phi_s(\cdot)$. It is proved in \cite{proudintegapp} that
$L_s=\frac{L}{d(D_s)}.$ We will consider instead of the function
$\Phi_s(\cdot)$  the function
$\tilde{\Phi}_{s}(\cdot):\mathbb{R}^n \rar \mathbb{R}$:
$$
\tilde{\Phi}_{s}(y,x)=\Phi_s(y)+  {L}_{s} \| y-x \|^2,
$$
for any fixed point $x \in \mathbb{R}^n$ and $y \in \mathbb{R}^n$.

As a result, the inequality \be {L_s} \| z \|^2 \leq (\grad^2
\tilde{\Phi}_{s}(x,x)z,z) \leq 3{L}_{s} \|  z \|^2 \,\,\,\,
\forall z \in \mathbb{R}^n, \label{newopt7a} \ee is true where
$\grad^2 \tilde{\Phi}_{s}(\cdot , x)= \tilde{\Phi''}_{s}(\cdot,x)$
is the matrix of the second mixed derivatives of the function
$\tilde{\Phi}_{s}(\cdot,x)$ with respect to the variable $y$.

Note that if the function $\Phi_s(\cdot)$  is bounded below, then
the function $\tilde{\Phi}_s(\cdot,x)$  is also bounded below for
any points $x$ and $y$ from $\mathbb{R}^n$. Also, it is clear that
$\grad \tilde{\Phi}_s(x,x) = \grad \Phi_s(x)$, where $\grad
\tilde{\Phi}_s(x,x)$ is the gradient of the function
$\tilde{\Phi}_s(\cdot,x)$ at the point $y = x$.

We assume that the functions $f(\cdot)$ and $\Phi_s(\cdot)$ are
bounded below and reach their infimum at some points.

{\bc \bf  Search method for a stationary point\ec}

Let the point $x_k$ at the $k$- th step have already been built.
Construct the point $x_{k+1}$. We put by definition
$\tilde{\Phi}_{s,k}(\cdot)= \tilde{\Phi}_s(\cdot, x_k)$.

1. Calculate $\lapl_k = -(\grad^2 \tilde{\Phi}_{s,k}(x_k))^{-1}
\grad \tilde{\Phi}_{s,k} (x_k).$

2. Find a non-negative integer $l_k$ for which \be
\tilde{\Phi}_{s,k}(x_k+2^{-l_k} \lapl_k) \leq \tilde
{\Phi}_{s,k}(x_k)-2^{-2l_k}\frac{L_s}{2}  \| \lapl _k \pr^2.
\label{newopt9} \ee

3. We assume  $x_{k+1} = x_k +2^{-l_k} \lapl_k$, $k=k+1$.

4. With increasing $k$ we decrease $d(D_s)$ such that the
inequality \be \frac{3 \| \lapl_k \|}{d(D_s)}< \eps_k
\label{newopt9a} \ee holds for some sequence $\{ \eps_k \}$, where
$\eps_k \rar +0$.  Go to the step 1.

Let us show that $\| \lapl_k \| \rar +0 $ as $k \rar \infty$ and
the number $l_k$ mentioned  in operation $1$ exists. Expand the
function $\tilde {\Phi}_{s,k}(\cdot)$ in a neighborhood of the
point $x_k$ in the Taylor series \be \tilde {\Phi}_{s,k}(x_k+\al
\lapl_k)=\tilde {\Phi}_{s,k}(x_k)+\al(\grad \tilde
{\Phi}_{s,k}(x_k), \lapl_k)+o_{s,k}(\| \al \lapl_k \|),
\label{newopt10} \ee where $o_{s,k}(\| \cdot \|)$ is an uniformly
infinitesimal function  in $k$.

As soon as $\lapl_k=-(\grad^2 \tilde \Phi_{s,k}(x_k))^{-1} \grad
\tilde \Phi_{s,k}(x_k)$, then $\grad \tilde \Phi_{s,k}(x_k)=
-\grad^2 \tilde \Phi_{s,k}(x_k) \lapl_k$. Consequently, $(\grad
\tilde \Phi_{s,k}(x_k), \lapl_k)=-(\grad^2 \tilde \Phi_{s,k}(x_k)
\lapl_k,\lapl_k)$. Therefore,  we can rewrite (\ref{newopt10}) in
the form \be \tilde \Phi_{s,k}(x_k+\al \lapl_k)=\tilde
\Phi_{s,k}(x_k)-\al(\grad^2 \tilde \Phi_{s,k}(x_k)\lapl_k,
\lapl_k)+o_{s,k}(\| \al \lapl_k \|). \label{newopt11} \ee As soon
as $ o_{s,k}(\| \cdot \|)$ is an uniformly infinitesimal function
with respect to $k$, then the inequality
$$
 o_{s,k}(\al \| \lapl_k \| ) \leq  \frac{\al \| \lapl_k \|}
 {N_s(\al \| \lapl_k \|)}
$$
is true for large $k$ where $ N_s(\al \| \lapl_k \|) \rar \infty$
as $\al \| \lapl_k\| \rar 0$.

From (\ref{newopt11}) we have
$$
\tilde \Phi_{s,k}(x_k+\al \lapl_k) \leq \tilde \Phi_{s,k}(x_k)-\al
{L_s} \| \lapl_k \|^2  +\frac{\al \| \lapl_k \|}{N_s(\al \|
\lapl_k \|)}=
$$
$$
= \Phi_{s,k}(x_k)-\al \| \lapl_k \|( {L_s}\| \lapl_k \|
-\frac{1}{N_s(\al \| \lapl_k \|)}).
$$
The value $\frac{1}{N_s(\al \| \lapl_k \|)}$ tends to zero as $\al
\| \lapl_k \| \rar 0$. Therefore, for small $ \|\lapl_k \|$ and,
consequently, for small $\| \grad \tilde \Phi_{s,k}(x_k) \|$, we
get \be {L_s}\| \lapl_k \| -\frac{1}{N_s(\al \| \lapl_k \|)} \meq
\frac{L_s \| \lapl_k \|}{2}. \label{newopt12} \ee

It follows from here that  the inequality \be \tilde
\Phi_{s,k}(x_k+\al \lapl_k) \leq \tilde \Phi_{s,k}(x_k)-\al
\frac{{L_s}}{2} \| \lapl_k \|^2 \label{newopt13} \ee is true for
sufficiently small $\| \grad \tilde \Phi_{s,k}(x_k) \|$ and any
$\al \in [0,1]$. Therefore, $\| \grad \tilde \Phi_{s,k}(x_k) \|$
tends to zero as $k \rar \infty$, since otherwise, as follows from
(\ref{newopt13}), the function $\tilde \Phi_{s,k} (\cdot)$ would
decrease in value $ \al \frac{{L_s}}{2} \| \lapl_k \|^2 $ along
the direction $\lapl_k$ at $k$ -th step. The last thing
contradicts to the lower boundedness of the function $\tilde
\Phi_{s,k}(\cdot)$ for all $k$ and $s$.

We will show that when the requirements of the step $4$ are
fulfilled, the function $ o_{s,k}(\cdot) $ is uniformly
infinitesimal in $k$ and $s$. From (\ref{newopt10}) for $\al =1$
we have \be o_{s,k}(\| \lapl_k \|) = \tilde
{\Phi}_{s,k}(x_k+\lapl_k) - \tilde {\Phi}_{s,k}(x_k)- (\grad
\tilde {\Phi}_{s,k}(x_k), \lapl_k). \label{newopt13a} \ee

We will use the midpoint theorem. Then
$$
\tilde {\Phi}_{s,k}(x_k+\lapl_k) - \tilde {\Phi}_{s,k}(x_k)= (
\grad \tilde {\Phi}_{s,k}(x_k+ \xi \lapl_k), \lapl_k)
$$
for $\xi \in [0,1]. $ Substitute the received expression in
(\ref{newopt13a}). We will have
$$
o_{s,k}(\| \lapl_k \|) = (\grad \tilde {\Phi}_{s,k}(x_k+ \xi
\lapl_k),\lapl_k) - (\grad \tilde {\Phi}_{s,k}(x_k), \lapl_k).
$$
We use the midpoint theorem again for the derivatives of the
function $\tilde {\Phi}_{s,k}(\cdot):$
$$
(\grad \tilde {\Phi}_{s,k}(x_k+ \xi \lapl_k),\lapl_k) - (\grad
\tilde {\Phi}_{s,k}(x_k), \lapl_k)=(\grad^2 \tilde
{\Phi}_{s,k}(x_k+ \eta \lapl_k)\lapl_k, \lapl_k).
$$
Therefore
$$
| o_{s,k}(\| \lapl_k \|) | \leq | \grad^2 \tilde{\Phi}_{s,k}(x_k+
\eta \lapl_k)\lapl_k, \lapl_k) |.
$$
It follows from the Lipschitz quality of the gradient $\grad
\Phi_{s}(\cdot) $ with the constant $L_{s}=\frac{L}{d(D_s)}$ that
the next evaluation
$$
\frac{| o_{s,k}(\| \lapl_k \|) |}{ \| \lapl_k \|} \leq \frac{3 \|
\lapl_k \|}{d(D_s)} < \eps_k.
$$
is true if (\ref{newopt7a}) is satisfied.

It follows from here that the functions $ o_{s,k}(\cdot) $  and
$\frac{1}{N_s(\al \| \lapl_k \|)}$ are uniformly infinitesimal
with respect to $k$ and $s$. Therefore, for small $\| \lapl_k \| $
the inequality (\ref{newopt12}) will be correct for $\al=1$.
Consequently, the inequality (\ref{newopt9}) is satisfied for
$l_k=0$ and the process goes with the full step $\lapl_k $.

\begin{thm} Any limit point of the sequence $\{ x_k \}$,
constructed according to the algorithm 1-4, is a stationary point
of the function $f(\cdot)$.
\end{thm}
{\bf Proof.} We have already proved that for small $\| \lapl_k \|
$ the process goes with the full step $\lapl_k.$ Since the
functions $\tilde \Phi_{s,k} (\cdot)$ are bounded below in
aggregate on $k, s$ and the inequality (\ref{newopt13}) is true
for all $k$ and $s$, then $\| \lapl_k \| \rar 0$  and $ \grad
\tilde{\Phi}_{s,k} (\cdot) \rar 0 $ for $s,k \rar \infty$.
Therefore, the sequence $\{ x_k \}$ has the limit points.

The following equalities
$$
\lapl_{k+1} = -(\grad^2 \tilde{\Phi}_{s,k+1}(x_{k+1}))^{-1} \grad
\tilde{\Phi}_{s,k+1}  (x_{k+1}), \,\,
\grad\tilde\Phi_{s,k}(x_{k+1})=o_{s,k}(\|\lapl_k\|),
$$
are correct where all $o_{s,k}(\cdot)$ in (\ref{newopt10}) are
uniformly infinitesimal in $k, s$ .

It follows from the definition of the function $ {\Phi}_{s}
(\cdot)$ that the gradient $ \grad {\Phi}_{s}  (\cdot)$ is a
convex hull of the generalized gradients of the function
$f(\cdot)$.

Taking into account what is said above about $\| \lapl_k \|$ and
$\grad \tilde{\Phi}_{s,k} (\cdot)$, and also from
uppersemicontinuity of the Clarke subdifferential mapping
\cite{demrub1}, \cite{clarke} we can imply that the inclusion $0
\in \p_{CL} f(x^*) $ is correct at a limit point $x_*$, i.e. $x_*$
is the stationary point of the function $f(\cdot)$. The theorem is
proved.$\Box$

To estimate the rate of convergence, we assume that $f(\cdot)$ is
convex and almost everywhere
$$
m \leq \| \grad^2 f(x) \| \leq M.
$$
From \cite{proudintegapp1} it follows that $\Phi_s(\cdot)$ is also
convex and for some $m_s, M_s >0$
$$
m_s \leq \| \grad^2 \Phi_s(x) \| \leq L_{s}.
$$
Define the function
$$
\tilde{\Phi}_{s,k} (y,x_k)=\Phi_s(y)+  {L}_s \| y-x_k \|^2
$$
for each $k$ and $y \in \mathbb{R}^n$ where $\tilde{L}_k > 0$ is
positive number depending on $k$ and  tending to zero as $k \rar
\infty$. To search for a stationary point of the function
$f(\cdot)$, we use the algorithm described below.

Let
$$
m_{s,k} \leq \| \grad^2 \tilde{\Phi}_{s,k} \| \leq L_{s,k}
$$
Since $m_{s,k} \rar m $ as $s(k) \rar \infty $, we assume that
$m_{s.k} \leq m $ for all $s(k)$.

We first introduce {\em the conditions of coherence}, which give
to us the rules of coherent striving to infinity of the parameters
$k$ and $s$. We will write them briefly in the form of dependence
$s = s ( k )$. Denote by $L_{s(k)}$ the constant bounding from
above the norm of the matrix $\grad^2
\tilde{\Phi}_{s(k),k}(\cdot): \| \grad^2
\tilde{\Phi}_{s(k),k}(\cdot) \| \leq L_{s,k}$. During the process
of optimization we satisfy to {\em conditions of coherence}:

\begin{enumerate}
\item $L_{s,k} \| \Delta _k \| \rar 0$ as $s(k) \rar \infty$;
\item for convergence with superlinear rate, we require that
$$
q_{s(k),k}=\frac{m^{-1}}{N_{s,k} ( \| \Delta _k \| )} \rar 0
$$
as $s(k),k \rar \infty$, where $ \frac{\|  \Delta_k \|}{N_{s,k} (
\| \Delta _k \| )}$ is a upper bound of the function
$o_{s(k),k}(\cdot),$ obtained from the expansion of the function
$\tilde{\Phi}_{s,k} (\cdot)$ at the $k$-th step (\ref{newopt10}).
It is clear that $N_{s(k),k}( \| \Delta _k \| ) \rar \infty$ as $k
\rar \infty$.
\end{enumerate}

The conditions 1 and 2 can be easy satisfied. At first the
optimization process goes on with constant $s$. As soon as the
step size $ \| \Delta_k \| $ becomes quite small, that means large
enough $N_{s(k),k} ( \| \Delta _k \| )$, we increase $s$, decrease
diameter $d(D_s)$ and, consequently, increase $L_{s(k)}$ so that
to satisfy to the conditions of coherence 1 and 2. As we shall see
below, $q_{s(k)}$ is the coefficient of proportionality between
$\| \Delta _{k+1} \|$ and $\| \Delta _{k} \|$. Therefore, we are
able to evaluate $q_{s(k)}$ by the  coefficient of proportionality
between $\| \Delta _{k+1} \|$ and $\| \Delta _{k} \|$ and,
therefore, to satisfy to the clause 2 of the consistency
conditions.

{\bc \bf Superlinear optimization method for finding the minimum
point of any final convex function $ f (\cdot) $ \ec}

Let a point $x_k$ already been found. Construct the pint
$x_{k+1}$.

1. Calculate the $k$-th step.
$$
\lapl_{k} = -(\grad^2 \tilde{\Phi}_{s,k}(x_{k}))^{-1} \grad
\tilde{\Phi}_{s,k}(x_{k}).
$$

2. Find a non-negative integer $l_k$ for which
$$
\tilde{\Phi}_{s,k}(x_k+2^{-l_k} \lapl_k) \leq \tilde
{\Phi}_{s,k}(x_k)-2^{-2l_k}\frac{m_{s,k}}{2}  \| \lapl _k \pr^2.
$$

3. We put $x_{k+1} = x_k +2^{-l_k} \lapl_k$, $k=k+1$.

4. Calculate for $k=k+1$
$$
\lapl_{k+1} = -(\grad^2 \tilde{\Phi}_{s,k+1}(x_{k+1}))^{-1} \grad
\tilde{\Phi}_{s,k+1}  (x_{k+1}).
$$

5. If
$$
L_{s} \| \Delta_{k+1} \| \leq \eps_{k+1}
$$
for an arbitrarily chosen sequence $\{ \eps_k \}, \eps_k \rar +0,$
then we increase $s$ such that the inequality
$$
\| \Delta_{k+1} \| \leq q_{s,k} \| \Delta_k \|
$$
remained in force.

6. Go to the step 1 and continue until the step size becomes less
than the specified value.

Let us prove that the sequence $\{ x_k \}$ converges to a minimum
point of the function $f(\cdot)$ with superlinear speed.

\begin{thm} The sequence $\{ x_k \}$, constructed according to the algorithm 1-3,
converges to an unique stationary point $x_*$ of the function
$\Phi(\cdot)$. For large $k$ the following estimate for the rate
of convergence of the method is correct \be \| x_k-x_* \| \leq
\nu^k(\lapl_k) \| x_1-x_*\|, \label{newopt14} \ee where
$\nu(\lapl_k)\rar_k0$ as $\| \lapl_k \| \rar_k0$.
\end{thm}

{\bf Proof}. As above, we are able to show that for sufficiently
large k  the process goes with a full step, i.e. $l_k=0$. From the
decomposition
$$
\grad \tilde{\Phi}_{s,k}(x_k+\lapl_k) = \grad
\tilde{\Phi}_{s,k}(x_k)+\grad^2 \tilde{\Phi}_{s,k}(x_k)\lapl_k
+o_{s,k}(\|\lapl_k\| )
$$
for
$$
\lapl_{k} = -(\grad^2 \tilde{\Phi}_{s,k}(x_{k}))^{-1} \grad
\tilde{\Phi}_{s,k}  (x_{k})
$$
we have
$$
\grad\tilde\Phi_{s,k}(x_{k+1})=o_{s,k}(\|\lapl_k\|).
$$
It is easy to check that
$$
\grad\tilde\Phi_{s,k}(x_{k+1})- \grad\Phi_s(x_{k+1})=\tilde
o_{s,k}(\|\lapl_k\|).
$$
But it is obvious that
$\grad\Phi_s(x_{k+1})=\grad\tilde\Phi_{s,k+1}(x_{k+1})$. Therefore
$\grad\tilde\Phi_{s,k+1}(x_{k+1})= \hat o_{s,k}(\|\lapl_k\|)$.
Since the function $\tilde{\Phi}_{s,k}(\cdot)$ has the continuous
second derivative, satisfying a Lipschitz condition, then
$o_{s,k}(\cdot), \tilde o_{s,k}(\cdot), \hat o_{s,k}(\cdot)$ are
the uniformly infinitesimal functions in $k$. From here
$$
\| \grad\tilde\Phi_{s,k+1}(x_{k+1}) \|  = \|\hat
o_{s,k}(\|\lapl_k\|)\| \leq \frac{ \| \lapl_k \| }{N_{s,k}(\|
\lapl_k \|)}.
$$
From the expression
$$
\lapl_{k+1} = -(\grad^2 \tilde{\Phi}_{s,k+1}(x_{k+1}))^{-1} \grad
\tilde{\Phi}_{s,k+1}  (x_{k+1}).
$$
we have the evaluation
$$
\| \lapl_{k+1} \| \leq \frac{m^{-1}}{N_{s,k}(\| \lapl_k \|)} \|
\lapl_k \|,
$$
where $ N_{s,k}(\| \lapl_k \|) \rar \infty$, as $\| \lapl_k \|
\rar 0$. For large $k$ we achieve that the inequality
$$
0< \frac{m^{-1}}{N_{s,k}(\| \lapl_k \|)}<1
$$
is correct (the condition of coherence). Therefore, the sequence
$\{ x_k \}$ converges to a single point $x_*$  and
$$
\| x_{k+1} - x_* \| \leq \sum_{i=k+1}^{\infty} \| \lapl_i \| =
\frac{(m^{-1}/ N(\|  \lapl_k \|)) \| \lapl_k \|}{1- m^{-1}/ N(\|
\lapl_k \|)}.
$$
As soon as
$$
\| \lapl_k \| \leq  (\frac{m^{-1}}{N_{s,k}(\| \lapl_k \|)})^k \|
\lapl_1 \|,
$$
then
$$
\| x_{k+1}-x_* \| \leq \frac{(m^{-1}/ N_{s,k}(\| \lapl_k
\|))^{k+1} \| \lapl_1 \|}{1-  {m^{-1}/ N_{s,k}(\| \lapl_k \|)}}.
$$

Thus, the inequality (\ref{newopt14}) is proved.  $\Box$

\begin{rem}\hspace{-2mm}.
The inequality (\ref{newopt14}) proves the superlinear convergence
rate of the optimization method. Indeed, the coefficient between
$\| x_{k+1}-x_* \|$ and $\| x_1-x_* \|$ is equal to $q_k^k$, where
$q_k \rar 0,$ as $k \rar \infty$.
\end{rem}


\section{Conclusion}

The methods for finding for a stationary point of Lipschitz
function and a minimum point of arbitrary convex function are
proposed in this paper. To achieve a high rate of convergence, it
is necessary to make consistent reduction  of the diameter $d(D)$
of the set $D$, which the integral averaging is doing on, with
decreasing the length of step of optimization process. Rules for
consistent reduction of the lengthes of steps and the diameters of
the sets $D_s$ are given.

\end{document}